\documentclass[12pt,a4paper]{article}

\usepackage[latin1]{inputenc}
\usepackage{amsmath}
\usepackage{amsfonts}
\usepackage{amssymb}
\usepackage{xspace}
\usepackage{graphics}

\renewcommand\iff{if and only if\xspace}
\renewcommand\th{\theta}
\def\l{\left}  \def\r{\right}

\newcommand\qed{\hfill$\quad${\rule{3mm}{3mm}}\medskip\\}
\newcommand\LB[1]{\label{#1}}
\newcommand\BE[2]{\begin{#1} #2 \end{#1}}
\newcommand\EQ[2]{\BE{equation}{\LB{#1} #2}}
\newcommand\ARR[2]{\BE{array}{{#1} #2}}
\newcommand\EQA[3]{\EQ{#1}{\ARR{#2}{#3}}}

\newcommand\ds{\displaystyle}

\def\a{\alpha}
\def\b{\beta}
\newcommand\f{\varphi}

\newcommand\Gm{\Gamma}

\newcommand \lm{\lambda}
\newcommand \Lm{\Lambda}

\newcommand\eps{\varepsilon}

\newcommand\dl{\delta}

\newcommand\lmn{\lambda_n}

\newcommand{\bbN}{\blackboard{N}}

\newcommand{\blackboard}[1]{\mathbb#1}

\newcommand{\iy}{\infty}

\newcommand\lmm{\lm_{nj}}

\title{%
{\bf Spectrum of the heat equation with memory\\
}
}

\author{
{S. Ivanov\thanks{Saint-Petersburg State University, Chemical department,
Ul'yanovskaya 1, Petrodvorets, St. Petersburg, 198904, Russia
{\tt sergei.ivanov@pobox.spbu.ru.}{Research of Sergei Ivanov
was supported in part by the Russia Foundation for Basic Research, grant 08-01-00595a.}}}\and
{T. Sheronova \thanks{State Marine Technical University of St. Petersburg
}}
}
\date{~}

\begin{document}

\maketitle

\begin{abstract}
We consider the system
\begin{equation*}
\partial_t\theta(x,t)=\int_0^t k(t-s) \partial_{xx}\theta(x,s)\,d s \  x\in(0,\pi), \ t>0
\qquad \theta(0,\cdot)=\xi(\cdot),
\end{equation*}
with homogeneous Dirichlet boundary condition.
Here
$$
k(t)=\sum_1^\iy  a_k e^{-b_k t}
$$
with positive $a_k$, $0\le b_1<b_2\dots $ and
$$
\sum_1^\iy  a_k<\iy, \ b_k\uparrow +\iy.
$$
Assuming an additional condition to $b_k$, e.g.,
$b_{k+1}-b_k\ge \dl>0$, we obtain the structure of the spectrum
of the system.
\end{abstract}




\vskip1cm

\section{\LB{intro_etc} Introduction. Notations. Main results}

Gurtin and Pipkin in \cite{GuPip} introduce a model of heat transfer
with finite propagation speed. A linearized model with the zero
memory at $t=0$ can be written \cite{P05} as:
\begin{equation}\label{1}
\dot \theta(x,t)=\int_0^t k(t-s) \theta''(x,s)\,d s, \  x\in(0,\pi), \ t>0,
\qquad \theta(0,\cdot)=\xi(\cdot)\in L^2(0,\pi).
\end{equation}
with homogeneous DBC. Here $\dot \theta=\frac\partial{\partial t} \th$,
$\theta'=\frac\partial{\partial x} \th$. Regularity of this equation in more
general setting is studied in \cite{P05}, but in one-dimensional case
it is possible to  obtain
regularity results directly using Fourier method. Two kind of controllability of
this system are
studied in \cite{P05}, \cite{IP}.

We will study this equation with the kernel in the form
$$
k(t)=\sum_1^\iy  a_k e^{-b_k t}
$$
with positive $a_k$, $0\le b_1<b_2\dots $ and
\EQ{ab}{
\sum_1^\iy  a_k=\a^2<\iy, \ b_k\uparrow +\iy.
}

Note that for $k(t)=c^2$ we have in fact the wave equation $\ddot\th=c^2\th''$.

First, apply the Fourier method:
we set $\f_n=\sqrt{\frac2\pi}\sin nx$ and expand the solution and the initial
data in series in $\f_n$
$$
\th(x,t)=\sum_1^\iy \th_n(t)\f_n(x),\qquad \xi(x)=\sum_1^\iy \xi_n \f_n(x).
$$
For the components we obtain
\begin{equation}\label{thn}
\dot \theta_n(t)=-n^2\int_0^t k(t-s) \theta_n (s)d s. \  t>0,
\qquad \theta_n(0)=\xi_n.
\end{equation}

We will denote the Laplace image by the capital characters.
Applying the Laplace Transform to  System \eqref{thn} we find
$$
z\Theta_n(z)-\xi_n=-n^2 K(z) \Theta_n(z)
$$
or
$$
\Theta_n(z)=\frac1{z+n^2K(z)}\xi_{n}.
$$

We will need to study zeros of the function $G_n(z)=z+n^2K(z)$.
Let $\Lm_n$, $n=1,2,\dots,$
be the set of zeros  of $G_n$, $\Lm:=\{\Lm_n\}$.

\BE{definition}{
The set $\Lm$  is called the spectrum of System \eqref1.
}

The Laplace image of $k(t)$ is
$$
K(z)=\sum_1^\iy \frac{a_k}{z+b_k}
$$
This is a meromorphic  function real valued on the real axis and mapping
the upper half plane to the lower one and visa versa.
\BE{lemma}
{
$K(z)$  has only real zeros ${-\mu_j}$, $j=1,2,\dots,$ such that
$$
b_1<\mu_1< b_2< \mu_2 <b_3<\dots .
$$
}
Proof of the lemma
Since
$$
\Im K(x+iy)=-y\sum_1^\iy \frac{a_k}{(x+b_k)^2+y^2},
$$
all zeros are real.
The lemma follows now from from the fact that $K(x)$ runs the real axis on every interval
$(-b_j,-b_{j-1})$, see the figure. \qed

\includegraphics{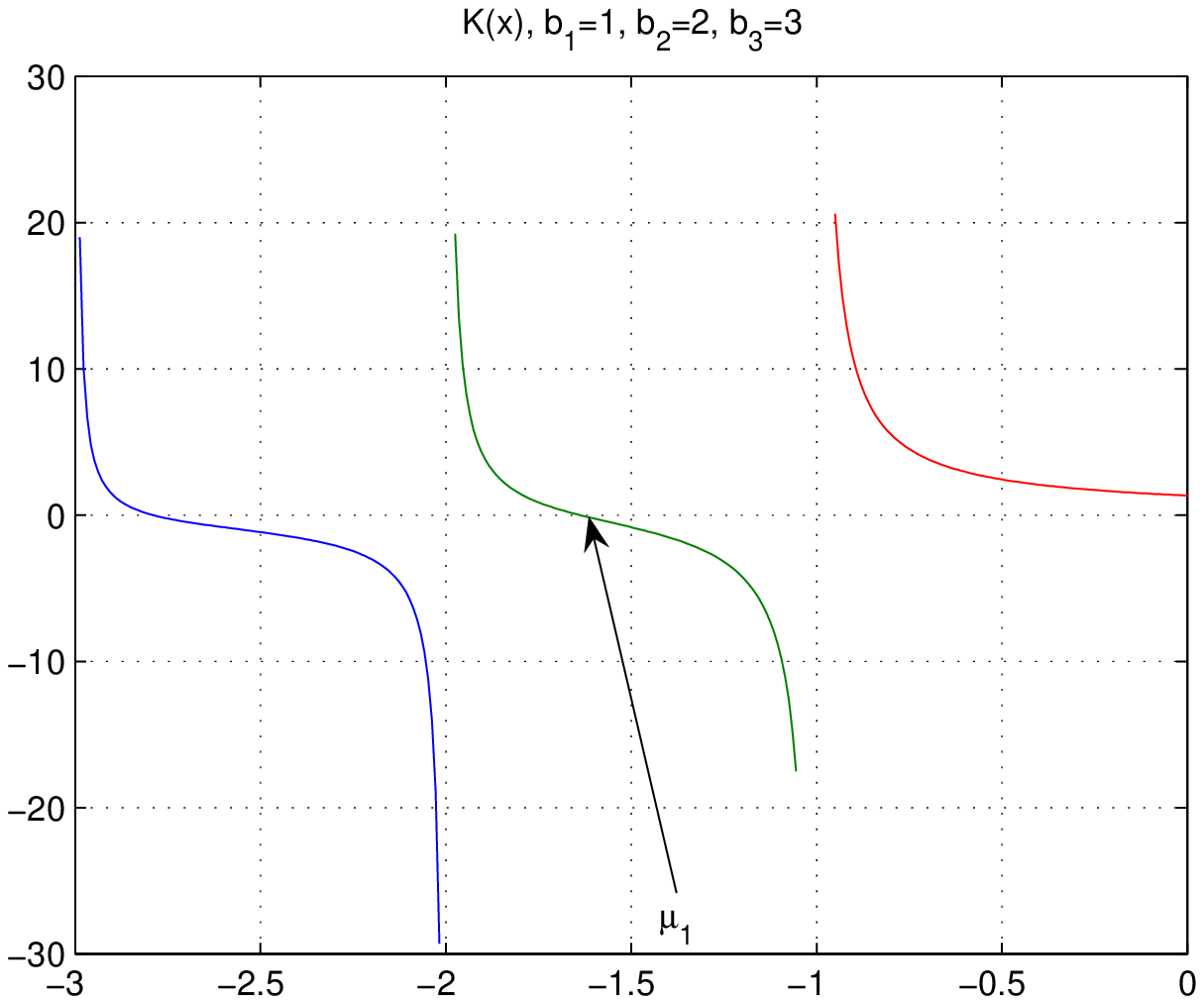}

\BE{theorem}{\LB{main}
Let
\EQ{dl}{
\sup_k\{b_k(b_{k+1}-b_k)\}=\iy.
}
Then the sets $\Lm_n$ can be represented as
$$
\Lm_n=\{\lmm\}_{j=1}^\iy\cup\{\lmn^+\}\cup\{\lmn^-\},
$$
with
$$
(i) \ \lm_{nj}= -\mu_j +o(1), \ j=1,2,\dots,; \ \lm_n^+=-i\a n+o(n),
\ \lm_n^-=\overline{\lm_n^+},
$$
and the spectrum is in the left half plane.

(ii) $\{\lmm\}$ are real and for fixed $j$ we have
$\lm_{nj}\uparrow - \mu_j$ and the sequence $\{\lmm\}_{n=1}^\iy$
is in $(-b_{j+1},-\mu_j)$:
$$
-b_{j+1}<\lm_{n,1}<\lm_{n,2}<\cdots<\lm_{n,k}<\cdots<-\mu_j.
$$
}
\BE{remark}{\LB{r1} In \cite{VVV} regularity of the
Gurtin-Pipkin equation  of the second order in time is studied
in more general situation.
In \cite{IV} the authors  consider this equation with model kernel
and find the asymptotic of complex zeros of the
the function
$$
z^2/n^2+1-K(z), n\in \bbN,
$$
where
$$
a_k=1/k^\a, \ b_k=k^\b, \ 0 < \a \le 1, \ \a+\b>1.
$$
In this case
$$
\sum_1^\iy c_k=\iy,\  \sum_1^\iy \frac {c_k}{a_k}<\iy.
$$
}

\BE{remark}{The authors are grateful to Prof. A.E. Eremenko for very helpful consultations.
}

\section{\LB{s1} Proof of the main theorem}

First show that  the spectrum is in the left half plane. If $z+n^2K(z)=0$, $z=x+iy$,  then
$$
x+iy +n^2\sum \frac{a_k(x+b_k)}{(x+b_k)^2+y^2}- iyn^2\sum \frac{a_k}{(x+b_k)^2+y^2}=0
$$
For the real part we have
$$
x+n^2\sum \frac{a_k(x+b_k)}{(x+b_k)^2+y^2}=0.
$$
If $x\ge 0 $, then this expression is positive.

We are going to apply the  Argument Principle. Fix $n$ and
take the rectangle contour $\Gm$ with
the vertices $(\pm X,\pm Y)$, $X,Y>0$. Let
$$
f(z)=K(z),\ g(z)=z/n^2.
$$
We show that we can take $X$,  and $ Y$ in such a way that
\EQ{fg}{
\Big| f|_\Gm  \Big| <  \Big| g|_\Gm  \Big|.
}
Consider the side $\Gm_1=\{\Re z=-X, |\Im z |\le Y\}$.
For $z\in \Gm_1$ we have
\EQA{N}{l}{
\ds
|K(z)|\le \sum\frac{a_k}{|-X+iy+b_k|}=
\le \sum\frac{a_k}{|-X+b_k|}=:q(X).
}

Take $X=X_N=(b_{N+1}+b_N)/2$ where we choose $N$ later.

\BE{lemma}{\LB{lm_dl}
\EQ{7_1}{
\frac1{X_N} q(X_N)\le \frac2{b_N\dl_N}\a^2.
}
}
Proof of the lemma.  From \eqref{N} we conclude
$$
q(X_N)=\sum_1^N \frac{a_k}{X_N-b_k}+ \sum_{N+1}^\iy  \frac{a_k}{b_k-X_N}.
$$
Since
$$
X_n-b_k=\frac12(b_{N+1}-b_k+b_N-b_k)\ge \frac12(b_{N+1}-b_k), \ k=1,2,\dots,N,
$$
and
$$
b_k-X_N\ge \frac12(b_k-b_N),\ k=N+1,N+2,\dots,
$$
we have
\EQ{KN}{
q(X_N)\le 2\l(\sum_1^N \frac{a_k}{b_{N+1}-b_k}+ \sum_{N+1}^\iy  \frac{a_k}{b_k-b_N}\r).
}
Set
$$
\dl_n=b_{n+1}-b_n.
$$
Evidently
$$
b_{N+1}-b_k\ge b_{N+1}-b_N= \dl_N, \, k=1,2,\dots,N,\ b_k-b_N\ge \dl_N, \, k=N+1,N+2,\dots,N.
$$
Now \eqref{KN} gives
$$
q(X_N)\le 2\frac1{\dl_N}\sum_1^\iy a_k=2\frac1{\dl_N}\a
$$
and then
$$
\frac1{X_N} q(X_N)\le \frac2{b_N\dl_N}\a^2.
$$
This proves the lemma.

Let us obtain \eqref{fg}  on $\Gm_1$.
Choose $N$ such that
\EQ{N1}{
\frac1{n^2}> \frac2{b_N\dl_N}\a^2,
}
what is possible by the assumption \eqref{dl}.
Then
$$
|g|_{\Gm_1}| \ge \frac {X_N}{n^2}\overset{\eqref{N1}}
> X_N  \frac2{b_N\dl_N}\a^2
\overset{\eqref{7_1}}
\ge
X_N \frac1{X_N} q(X_N)\ge |K|_{\Gm_1}|.
$$

Consider the side $\Gm_2=\{ |\Re z|<X, \ \Im z=Y \}$.
Here $\inf|g|=\inf |z/n^2|=Y/n^2$ and
\EQ{8_1}{
|K(z)|\le \sum\frac{a_k}{|x+iY+b_k|}\le \sum\frac{a_k}Y= \frac{\a^2} Y.
}
Let $Y>n\a$. Then
$$
\frac Y{n^2} > \frac{\a^2}{Y}$$
and \eqref{8_1} gives \eqref{fg}.

The estimate on $\overline{\Gm_2} $ is the same  as for $\Gm_2$ because
\EQ{sym}{
\bar K(z)=K(\bar z).
}
On the rest side $\Re z=X, \ |\Im z\le Y|$ we have
$$
|K|\le \sum \frac{a_k}X=\frac{\a^2} X.
$$
Since $|z/n^2|\ge X/n^2$ we have \eqref{fg} if
$X>n\a$.
Therefore, for $X,Y> n\a$ and with \eqref{N1}
we have \eqref{fg} and we can conclude that
$$
N(g)-P(g)=N(f+g)-P(f+g)
$$
where  $N$ and $P$ denote respectively the number of zeros and poles
of the function inside the contour $\Gm$,
with each zero and pole counted as many times as its multiplicity and order respectively.
Evidently, $N(g)-P(g)=1$. Inside $\Gm$ the function $f+g$ has $N$ simple poles at
$-b_N$, $-b_{N-1},\dots,-b_1$. Therefore $N(f+g)=N+1$.

Show that $G_n$ has $N-1$ real zeros $\lm_{nj}$
inside $\Gm_n$ satisfying
the theorem. For fixed $n$ the graph of one branch of $K(x)$ and of
$-x/n^2$ is

\includegraphics{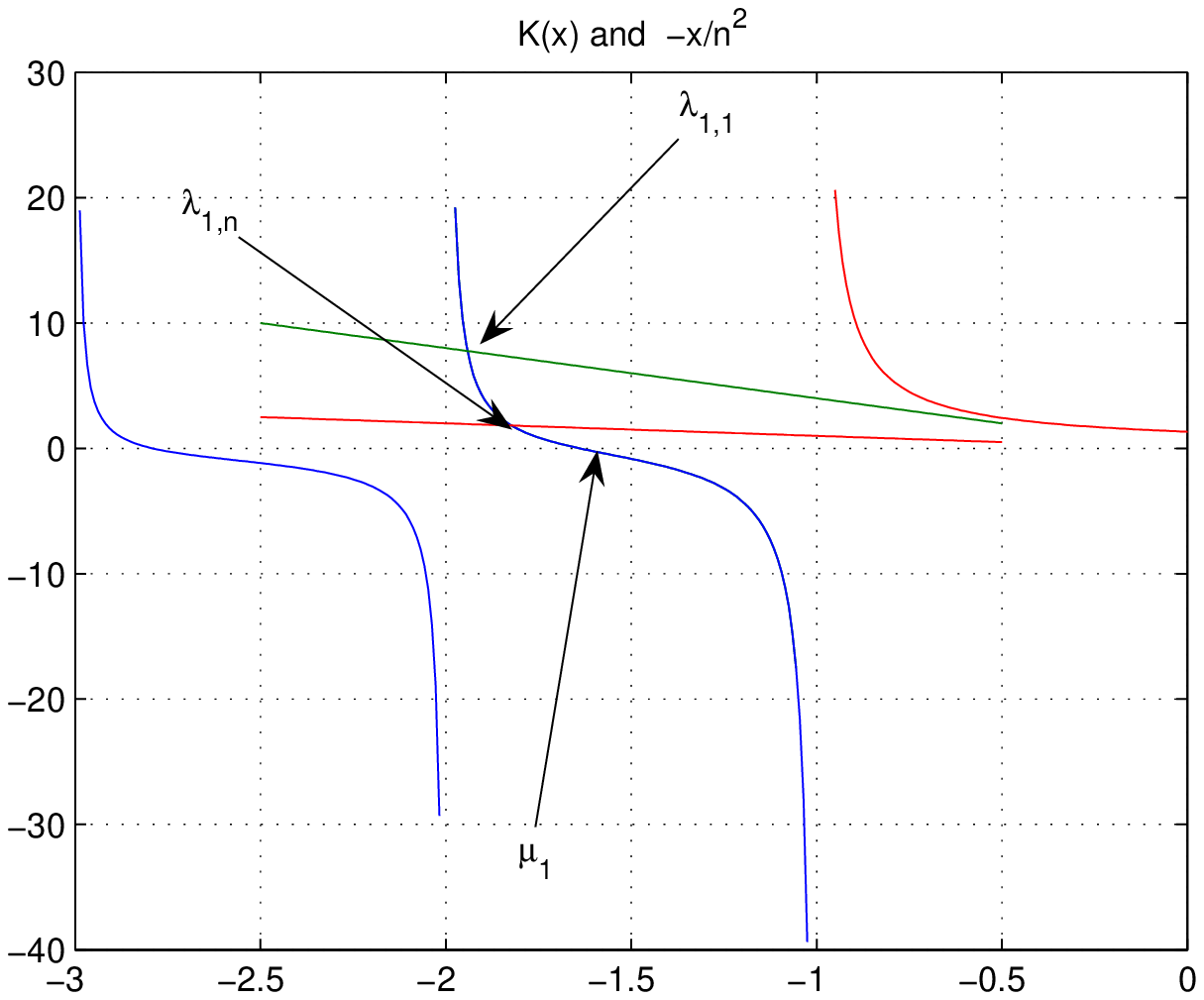}

We see that the straight line $-x/n^2$ intersects the graph of $K(x)$ in $(-b_{j+1},-\mu_j)$
and monotonically approaches $-\mu_j$ as $n\to \iy$, $j=1,2,\dots$.
Thus, there exists $N$ or $N-1$ real zeros of $G_n$ inside $\Gm_n$ depending
whether
$-(b_{N+2}+b_N)/2$ is more or less $\lm_{N,n}$. In the first case we have exactly
one complex zero of $G_n$ what is impossible by \eqref{sym}.

Thus, there exists two
(complex conjugated) zeros. Denote these zeros by $\lmn^+$  and $\lmn^-=\overline{\lmn^+}$.

Let us localize these roots. Consider the case $b_1=0$. Then for $z=iy$
$$
G_n(iy)=iy+n^2\sum_2^\iy \frac{a_kb_k}{b_k^2+y^2} - in^2 y\sum_1^\iy \frac{a_k}{b_k^2+y^2}.
$$
We find that $G_n(iy)=0$ only if $b_2=b_3=\dots=0$, i.e., $K(z)=\a^2/z$ and
$\lm_n^\pm=\pm i\sqrt{\a}n$. Consider the case $K(z) \ne \a^2/z$.

Take $\eps>0$ and the rectangle contour $G$ with
the vertices $(-\eps \a  n,in\a (1-\eps) )$,
$(-\eps\a  n,in\a (1+\eps) )$, $(\eps \a n,in\a (1+\eps) )$,
 $(\eps \a n,in\a (1-\eps) )$. Write $z/n^2+K(z)$ as
 $$
 \l[\frac z{n^2}  +\frac{\a^2} z\r]+\l[K(z) -\frac{\a ^2}z\r]=g+f.
 $$
 \BE{lemma}{\LB{asympK}
Let
\EQ{dl1}{
|\arg z|<\pi-\dl,
}
Then
\EQ{Kz}{
K(z)=\frac{\a^2} z(1+o(z)).
}

}
Proof of the lemma.
Show first that  for
\EQ{sum}{
|z+x^b|^2\asymp |z|^2+x^{2b},  \ x\ge0,\  b>0 .
}
Set $z=\rho e^{i\f}$, $\rho=|z|$
$$
|z+x^b|^2=(\rho\cos\f+x^b)^2+\rho^2\sin^2\f=\rho^2+2x^b\rho\cos\f+x^{2b}.
$$
Evidently
$$
\rho^2+2x^b\rho\cos\f+x^{2b}
\le
(\rho+x^{b})^2\le 2(\rho^2+x^{2b})
$$
and
$$
\rho^2+2x^b\rho\cos\f+x^{2b}
\ge
\rho^2-2x^b\rho\cos\dl+x^{2b}\ge
(\rho^2+x^{2b})-(\rho^2+x^{2b})\cos\dl
$$
$$
=
(1-\cos\dl)(\rho^2+x^{2b}).
$$
This gives \eqref{sum}

Now we can apply the Weierstrass theorem to the series
$$
zK(z)=\sum\frac{a_k z}{z+b_k}.
$$
Indeed, by \eqref{sum}
$$
\l|\frac{a_k z}{z+b_k}\r|\prec \frac{a_k|z|}{|z|+b_k}\le a_k
$$
and
$$
\frac{a_k z}{z+b_k}\ \to a_k, \ z\to\iy.
$$
Therefore
$$
zK(z)=\lim_{z\to\iy}\sum\frac{a_k z}{z+b_k}=\sum_1^\iy a_k=\a^2.
$$
The lemma is proved.

Find the estimate of $g|_G$ from below. For $\eps<1$
$$
|g|_G|=\frac{|z+i\a n||z-i\a n|}{|z|n^2}
\ge \frac{  n\a(2-\eps)\eps\a n}{|z||n^2}>\frac{\a^2 \eps}{|z|}.
$$
From Lemma \ref{asympK} we have in the sector \eqref{dl1}
$$
f(z)=h(z)\frac1z, \ h(z)=o(1).
$$
Take $n$ large enough in order that
$$
|h(z)|_G|<\eps.
$$
Then
$$
|g|_G|> |f|_G|.
$$
Now $N(f+g)=N(g)=1$. We have proved that $\lm_n^+$ is in the box centered
at $in\a$ and with the diameter  $O(\eps n)$.

The theorem is proved.
\BE{remark}{
If we know asymptotic or estimate of the parameters
$\{a_k,b_k\}$, we can find a more precise expression of
points of $\Lm$, especially of $\lm_n^\pm$, see \cite{IV}.
}
\BE{remark}{
The condition \eqref{dl} is evidently fulfilled
if $b_{k+1}-b_k\ge \dl>0$. For a model case $b_k=ck^\a$
$$
b_k(b_{k+1}-b_k)\asymp k^{2\a-1}
$$
and \eqref{dl} is true \iff $\a>1/2$.
Using the proposed approach we have to show that
$$
\frac1{b_N}\sum\frac{a_k}{|b_k-(b_{N+1}-b_N)/2|} \ \to \ 0, \ as \ N\to\iy.
$$
This is not true for slowly increasing $b_k$, say, $b_k=\log\log k$.
Nevertheless we have the conjecture that main theorem is true for every kernel satisfying
\eqref{ab}.
}

\end{document}